\documentclass[reqno,11pt]{amsart}
\usepackage{amsmath,amssymb,latexsym,soul,cite,mathrsfs}
\usepackage{color,enumitem,graphicx}
\usepackage[colorlinks=true,urlcolor=blue,
citecolor=red,linkcolor=blue,linktocpage,pdfpagelabels,
bookmarksnumbered,bookmarksopen]{hyperref}
\usepackage[english]{babel}
\usepackage[left=2.6cm,right=2.6cm,top=2.9cm,bottom=2.9cm]{geometry}
\usepackage{pgfplots}
\pgfplotsset{compat=1.15}
\usepackage{mathrsfs}
\usetikzlibrary{arrows}
\usepackage{geometry}
\usepackage{multirow}

\usepackage{float}
\usepackage{amssymb}
\usepackage{verbatim} 
\usepackage{amsmath}
\usepackage{amsthm}
\usepackage{epstopdf}
\usepackage{epsfig}
\usepackage[all]{xy}
\usepackage{enumerate}
\usepackage{hyperref}

\newtheorem{theorem}{Theorem}

\newtheorem{definition}[theorem]{Definition}

\newtheorem*{theorem*}{Theorem}
\newtheorem*{lemma*}{Lemma}
\newtheorem*{corollary*}{Corollary}
\newtheorem*{proposition*}{Proposition}
\newtheorem*{remark*}{Remark}
\newtheorem*{definition*}{Definition}
\newtheorem*{conjecture*}{Conjecture}
\newtheorem*{claim*}{Claim}

\DeclareMathOperator{\vol}{\rm vol}
\DeclareMathOperator{\area}{\rm area}


\begin{document}
\title[]{Area-minimizing unit vector fields on some spherical annuli}
\author{Fabiano Brito$^1$}
\author{Jackeline Conrado$^2$}
\author{Jo\~ao Lucas $^3$}
\author{Giovanni Nunes$^4$}

\address{Centro de Matem\'atica, Computa\c c\~{a}o e Cogni\c c\~{a}o, 
Universidade Federal do ABC, Santo Andr\'{e} - SP, 09210-170, Brazil}
\email{fabiano.brito@ufabc.edu.br}
\email{francisco.joao@ufabc.edu.br}

\address{Dpto. de Geometria e Representação Gráfica, Instituto de Matem\'{a}tica e Estat\'{i}stica, 
Universidade do Estado do Rio de Janeiro, Rua São Francisco Xavier, 524 - Pavilhão Reitor João Lyra Filho, Rio de Janeiro - RJ, 20550-900, Brazil}
\email{jackeline.conrado@ime.uerj.br}

\address{Dpto de Matem\'atica e Estat\'istica, Instituto de Matem\'{a}tica e F\'isica,
Universidade Federal de Pelotas, Rua Gomes Carneiro 1, Centro - Pelotas, 96010-610, Brazil}
\email{giovanni.nunes@ufpel.edu.br}



\begin{abstract}
We establish in this paper a sharp lower bound for the area of a unit vector field $V$ defined on some spherical annuli in the Euclidean sphere $\mathbb{S}^2$.
\end{abstract}

\maketitle

\section{Introduction }

The Sasaki metric, introduced by Shigeo Sasaki in [\textcolor{red}{7}, $1958$], provides a natural metric on the tangent bundle $TM$ of a Riemannian manifold $M$. This metric allows associating a volume with each unit vector field on $M$, considering them as sections of the tangent bundle $TM$ and regarding this volume as the measure induced by the Sasaki metric on these sections, viewed as submanifolds of $TM$. If $M$ is a surface and  $V$ a unit vector field on $M$. The \textit{\textbf{area of a unit vector field $V$}} on $M$ is defined as the area of the surface $V(M)$ in $T^1M$. 
For a comprehensive understanding of the geometry of the tangent bundle, we suggest consulting chapter $1$ of the book \textit{Harmonic Vector Fields: Variational Principles and Differential Geometry} by Sorin Dragomir and Domenico Perrone, [\textcolor{red}{5}, $2012$].

Throughout the years, results concerning the minimization of volume or area of unit vector fields on  $M$ have been obtained. In the next section, we shall present these results related to the area minimization of unit vector fields on $\mathbb{S}^2\backslash\{N\}$ and $\mathbb{S}^2\backslash\{N,S\}$.

In this work, we establish a minimization result for unit vector fields defined on spherical annuli in the Euclidean sphere $\mathbb{S}^2$. More precisely, we consider regions of the unit Euclidean sphere \( \mathbb{S}^2 \), bounded by two parallels circles. We consider unit vector fields \( V \) that satisfy the following properties:

    \noindent \textit{(i) The vector field \( V \) is tangent to the two parallel circles, which are oriented in opposite directions, with latitudes given by \( \alpha_0 \) and \( -\alpha_0 \), respectively.} \\
    \textit{(ii) The angle \( \theta \) that \( V \) makes with the parallels is constant along them.}\\
    \textit{(iii) The unit vector field $V$ is perpendicular to the meridian $S_0$.}\\

Considering these conditions, the main questions we explore are
\begin{itemize}
    \item [1.] To determine the lower bounds for the area of \( V \) on spherical annuli, where \( \alpha_0 \in (0, \pi/2) \).
    \item [2.] To identify and explicitly describe the unit vector fields that achieve this minimum and to provide a geometric and qualitative description of them for each choice of \( \alpha_0 \).
  
\end{itemize}
A precise definition of these vector fields on spherical annuli in the Euclidean sphere $\mathbb{S}^2$ is provided in the definition \ref{Def:AreaMinimizing_Annulus}. Therefore, we prove the following theorem
\begin{theorem}
Let $V$ be a unit vector field as described in definition \ref{Def:AreaMinimizing_Annulus}. Let $\theta(p)$ be the angle function of $V$ with respect to the parallels at the point $p=(\alpha, \beta)$ in $A_{\alpha_0}$. If $\theta(-\alpha_0, \beta) = 0$ and $\theta(\alpha_0, \beta) = \pi$, for all $\beta \in (0,2\pi]$. Then,
\begin{itemize}
    \item[(i)]  $area(V) \geq 2\pi \cos \left(\alpha_0\right) + K_{\alpha_0}$, where $K_{\alpha_0}$ is constant depending just on $\alpha_0 \in \left(0, \dfrac{\pi}{2}\right)$.\\
    \item[(ii)] If the angle function $\theta: [-\alpha_0, \alpha_0] \longrightarrow [0, +\infty)$ is given by
    \[
    \theta(\alpha) = \arcsin\bigg(\cot (\alpha_0) 
    \tan (\alpha)\bigg) + \frac{\pi}{2},
    \]
    then, the vector field 
    \[
    V(\alpha, \beta) = \cos \big(\theta(\alpha)\big) \, e_1 + \sin \big( \theta(\alpha)\big) \, e_2
    \]
    achieves the lower bound of (i).
\end{itemize}
    
\end{theorem}

\section{Area-minimizing unit vector fields on $\mathbb{S}^2\backslash\{N\}$ or $\mathbb{S}^2\backslash\{N,S\}$}


Let $M$ be a closed oriented Riemannian manifold and $V$ a unit vector field on $M$. Consider the unit tangent bundle $T^1M$ equipped with the Sasaki metric. The \textit{\textbf{volume of a unit vector field $V$}} is defined (see [\textcolor{red}{6}, $1986$]) as the volume of the submanifold $V(M)$, the image of the immersion $V: M \rightarrow T^1M$, 
\begin{equation*}\label{Eq:volume}
\vol(V):=\vol(V(M)).
\end{equation*}


\noindent If $M$ is a surface and  $V$ a unit vector field on $M$. The \textit{\textbf{area of a unit vector field $V$}} on $M$ is defined as the area of the surface $V(M)$ in $T^1M$, and it coincides with the volume of a unit vector field $V$
\begin{align*}
    \area(V) := \area(V(M)) = \vol(V).
\end{align*}

\noindent It is well known that there is no globally defined vector field on $\mathbb{S}^2$, so the infimum of the area is achieved only when there is at least one singularity. The Pontryagin vector field is a unit vector field with one singularity obtained by parallel translating a given vector along any great circle passing through a given point. In [\textcolor{red}{1}, $2010$], Vincent Borrelli and Olga Gil-Medrano showed that Pontryagin vector fields of the unit 2-sphere with one singularity are area-minimizing. 

In [\textcolor{red}{2}, $2008$] Fabiano Brito, Pablo Chac\'on and David L. Johnson established an explicit relationship between the volume of a unit vector field and its Poincaré index around isolated singularities. 
\begin{theorem}[Brito, Chacón and Johnson, \textcolor{red}{2}]\label{BCJ}
Let $M = \mathbb{S}^{n} \backslash \left\{N,S\right\}$, $n = 2$ or $3$, be the standard Euclidean sphere with two antipodal points $N$ and $S$ removed. Let $V$ be a unit vector field defined on $M$. Then,
\begin{itemize}
\item [i)] for $n=2$, $\vol(V)\geq \frac{1}{2} (\pi + |I_{V}(N)| + |I_{V}(S)| -2)\vol(\mathbb{S}^2)$,

\item [ii)]for $n=3$, $\vol(V)\geq (|I_{V}(N)| + |I_{V}(S)|)\vol(\mathbb{S}^3)$,
\end{itemize}
where $I_V(p)$ stands the Poincaré index of $V$ around $p$.
\end{theorem}

\noindent
Theorem~\ref{BCJ} has been extended to odd dimensional spheres $\mathbb{S}^{2n+1}$ as follows:

\begin{theorem}[Brito, Gomes and Gonçalves, \textcolor{red}{4}] If $V$ is a unit vector field on $\mathbb{S}^{2n+1}\backslash\{\pm p\}$, then
\begin{equation*}
    \vol(V)\geq \frac{\pi}{4}\vol(\mathbb{S}^{2n})(|I_{ V}(p)| + |I_{ V}(-p)|).
\end{equation*}
\end{theorem}

In [\textcolor{red}{3}, $2021$], Fabiano Brito et al., established a sharp lower bounds for the total area of unit vector fields on antipodally punctured unit Euclidean sphere, and these values depend on the indices of the singularities. 

\begin{theorem}[Brito \textit{et al.}, \textcolor{red}{3}]\label{Thm:BCGN}
Let $V$ be a unit vector field defined on $M=\mathbb{S}^2\backslash\{N, S\}$. If $k = \max\{I_V(N), I_V(S)\}$, $k\neq0$, $k\neq2$ then 
\begin{eqnarray*}
    \area(V) \geq \pi L(\epsilon_k)
\end{eqnarray*} where $L(\epsilon_k)$ denote the length of the ellipse $\dfrac{x^2}{k^2}+\dfrac{y^2}{(k-2)^2}=1$ and $I_V(p)$ stands for the Poincaré index of $V$ around $p$.
   
\end{theorem}

\noindent Furthermore, in [\textcolor{red}{3}, $2021$], the authors exhibited a family of the unit vector fields $V_k$ on  $\mathbb{S}^2\backslash\{N,S\}$ that are area-minimizing within each index class $k$. This family has the following characteristics:

\begin{itemize}
    \item [1.] $V_k$ turns $k-1$ times along each parallel at a constant angle speed with respect to the oriented orthonormal frame $\left\{e_1, e_2 \right\}$ on  $\mathbb{S}^2\backslash\{N,S\}$;
    \item [2.] $V_k$ is parallel along meridians, which means that
    \begin{equation*}
d\theta(e_2) = 0.
\end{equation*}
    \item [3.] $V_k$ makes an angle \( \theta \) with each parallel $\alpha \in \left( -\frac{\pi}{2}, \frac{\pi}{2} \right)$, and $\theta$ has constant variation along them, implying that 
\begin{equation*}\label{Eq:variacao_constante}
d\theta(e_1) = \frac{k-1}{\cos \alpha}.
\end{equation*}    
\end{itemize}

\begin{figure}[H]
		\centering
		\includegraphics[height=7cm]{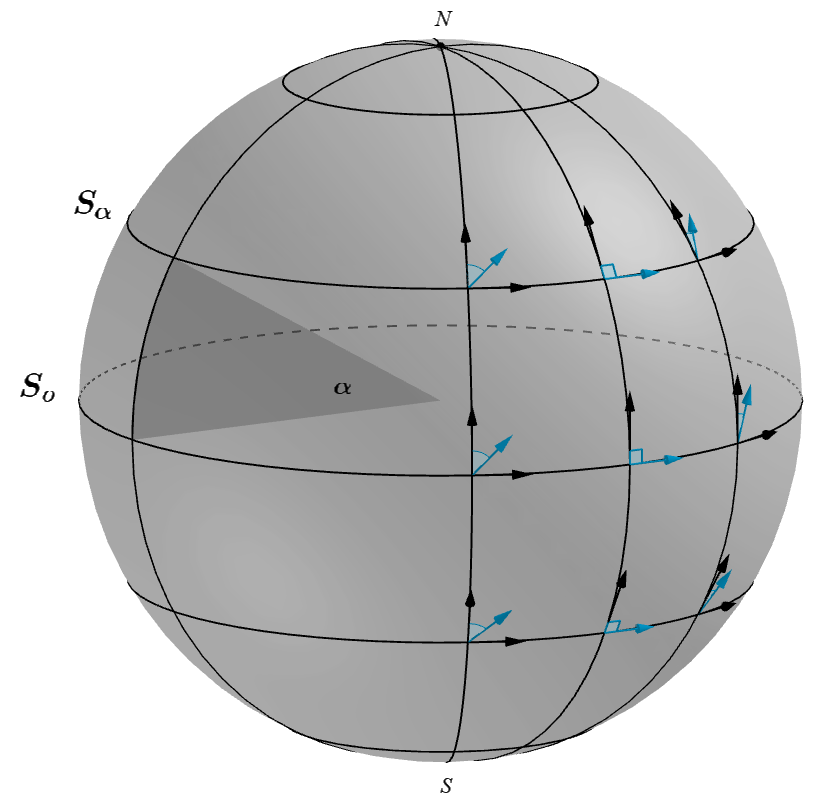}
		\vspace{0.5cm}
		\caption{Visual representation of the area-minimizing unit vector field $V_{k}$ on $\mathbb{S}^2\backslash\{N,S\}$}
		\label{Fig:Area_Minimizing}
	\end{figure}

In Figure \ref{Fig:Area_Minimizing} is a visual representation about the behavior of area-minimizing vector field $V_k$ on $\mathbb{S}^2\backslash\{N,S\}$ (in blue), where $S_{\alpha}$ and $S_0$ denote the parallel $\alpha$ and the Equator, respectively. 

Consider another oriented orthonormal local frame $\{V_k, V_k^\perp \}$ on $\mathbb{S}^2\backslash\left\{ N,S\right\}$ compatible with the orientation of $\{e_1, e_2\}$. The area of $V_k$ is given by
\begin{eqnarray}
    \area \left(V_k\right) = \int_{\mathbb{S}^2} \sqrt{1+ \gamma^2 + \delta^2}\,\,\nu,
\end{eqnarray}
where $\nu$ is the volume form, $\gamma = g\left(\nabla_{_{V_k}} V_k \,, V_k^\perp\right)$ and $\delta = g\left(\nabla_{_{V_k^\perp}} V_k^\perp \, , V_k\right)$ are the geodesic curvatures associated to $V_k$ and $V_k^\perp$, respectively.

\noindent Let $\theta \in [0, 2\pi)$ be the oriented angle from $e_2$ to $V_k$ and $\alpha \in \left[ -\frac{\pi}{2}, \frac{\pi}{2} \right]$. If $V_k = \cos\left(\theta\right)e_1 + \sin\left(\theta\right)e_2$ and $V_k^\perp = -\sin\left(\theta\right)e_1 + \cos\left(\theta\right)e_2$, then
\begin{equation}\label{curva_geode_tangente}
\gamma = \cos \left(\theta\right) \big(\tan \alpha + d\theta(e_1)\big) + \sin\left(\theta\right) d\theta(e_2), 
\end{equation}
\begin{equation}\label{curva_geode_ortogonal}
\delta = \sin (\theta) \big(\tan \alpha + d\theta(e_1)\big) - \cos (\theta) d\theta(e_2),
\end{equation}
and 
\begin{eqnarray}\label{Eq:Integrando_Volume}
    1+\gamma^2+\delta^2 = 1 + \big( \tan \alpha + d\theta(e_1)\big)^2 + d\theta(e_2)^2.
\end{eqnarray}

\noindent The Lemma in [\textcolor{red}{3}] provides a straightforward computation from \eqref{curva_geode_tangente} to \eqref{Eq:Integrando_Volume}.

\noindent 
 The equation \eqref{Eq:Integrando_Volume} allows us to rewrite the area functional as an integral depending on the latitude $\alpha$ and the derivatives of $\theta$:
\begin{eqnarray}\label{Eq: Vol_derivatives}
    \area \left(V_k\right) = \int_{\mathbb{S}^2} \sqrt{1+  \big( \tan \alpha + d\theta(e_1)\big)^2 + d\theta(e_2)^2}\,\nu,
\end{eqnarray}

The formula for the area of the unit vector field obtained in equation \eqref{Eq: Vol_derivatives} will be used to demonstrate our area minimization result on the spherical annulus of $\mathbb{S}^2$.



\section{Proof of Theorem}

\noindent 
Let $\mathbb{S}^{2}\backslash\{N,S\}$ be the Euclidean sphere with two antipodal points, $N$ and $S$, removed. Denote by $g$ the usual metric of $\mathbb{S}^2$ induced from $\mathbb{R}^3$ and the oriented orthonormal frame $\left\{ e_1, e_2 \right\}$ on $\mathbb{S}^2\backslash\left\{ N,S\right\}$ such that $e_1$ is tangent to the parallels and $e_2$ to the meridians. Consider $S_\alpha$ be the parallel of $\mathbb{S}^2$ at latitude $\alpha \in \left( -\frac{\pi}{2}, \frac{\pi}{2} \right)$ and $S_\beta$ be the meridian of $\mathbb{S}^2$ at longitude $\beta \in (0, 2\pi]$. The \textbf{\textit{spherical annuli}} in $\mathbb{S}^2$ bounded by the parallels $S_{-\alpha_0}$ and $S_{\alpha_0}$ are defined as  
\begin{eqnarray*}
A\left(-\alpha_0, \alpha_0\right) = \{ (x, y, z) \in \mathbb{S}^2 \mid -\alpha_0 \leq \arcsin{(z)} \leq \alpha_0 \},     
\end{eqnarray*}
where $\alpha_0 \in \left( 0 , \dfrac{\pi}{2} \right)$ determines the latitudinal bounds of the annular region on the sphere.
Furthermore, since every point $p=(x,y,z)$ on $\mathbb{S}^2$ can be defined by the latitude $\alpha$ and the longitude $\beta$, in order to simplify the notation, we rewrite the \textbf{\textit{spherical annulus}} $A\left(-\alpha_0, \alpha_0\right)$ in $\mathbb{S}^2$ as follows
\begin{eqnarray}\label{Def:Annulus_region}
A_{\alpha_0} = \{ (\alpha, \beta) \in S^2 \mid -\alpha_0 \leq \alpha \leq \alpha_0 \},    
\end{eqnarray}
where $\alpha_0$ is a fixed constant in the interval $\left( 0 , \frac{\pi}{2} \right)$. 

\begin{definition}\label{Def:AreaMinimizing_Annulus}
Let $V$ be a unit vector field on the annulus region
$A_{\alpha_0}$, with $\alpha_0 \in \left( 0 , \frac{\pi}{2}\right)$ that satisfies the following:
\begin{itemize}
    \item [(i)] $V$ is tangent to the boundaries $S_{-\alpha_0}$ and $ S_{\alpha_0}$ of the spherical annulus $A_{\alpha_0}$, such that
    \begin{eqnarray*}
        V(-\alpha_0, \beta) = -V(\alpha_0, \beta),
    \end{eqnarray*}
            for all $\beta$ in $(0,2\pi]$.
    \item [(ii)] For all $\beta$ in $(0,2\pi]$, $V(0,\beta)$ is perpendicular to the parallel $S_0$.
    \item [(iii)] Let $\theta: A_{\alpha_0} \rightarrow \mathbb{R}$ be the angle function of $V$ with the parallels. The angle $\theta$ that $V$ makes with the parallels is constant along them, i.e.,
    \begin{eqnarray}\label{Eq:Angulo_cste_paralelos}
        d\theta_p\left(e_1\right) = \theta_1 (p) \equiv 0,
    \end{eqnarray}
    for all $p$ in $A_{\alpha_0}$, where $\left\{ e_1, e_2 \right\}$ is the oriented orthonormal frame aforementioned.
\end{itemize}
\end{definition}

\begin{figure}[H]
		\centering
		\includegraphics[height=7cm]{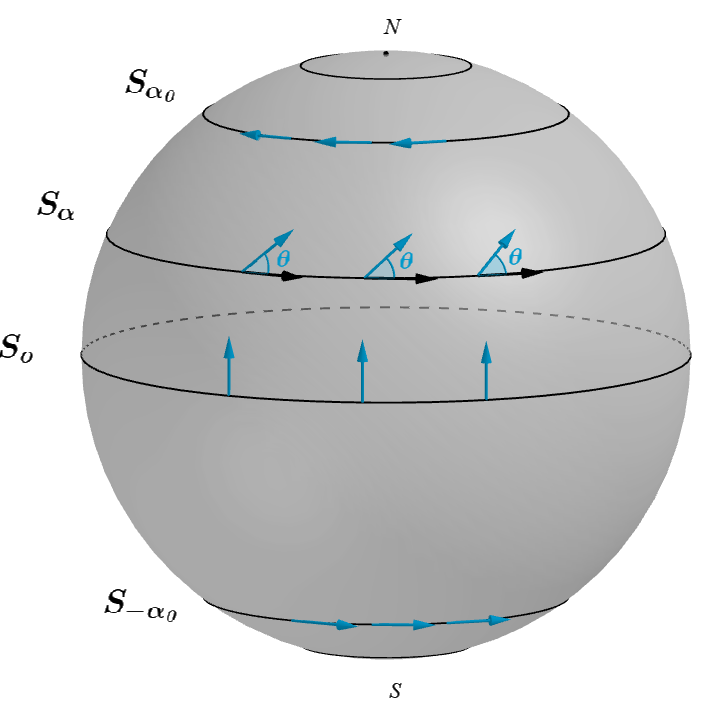}
		\vspace{0.5cm}
		\caption{Visual representation of the area-minimizing unit vector field $V$ on $A_{\alpha_0}$}
		\label{Fig:Area_Minimizing_aneis}
	\end{figure}

A visual representation about the behavior of area-minimizing vector field $V$ on $A_{\alpha_0}$ (in blue) is given by Figure \ref{Fig:Area_Minimizing_aneis} , where $S_{-\alpha_0}$ and $ S_{\alpha_0}$ are the boundaries of the spherical annulus $A_{\alpha_0}$,  $S_{\alpha}$ and $S_0$ denote the parallel $\alpha$ and the Equator, respectively, 

    \begin{proof} \textit{of Theorem $1$.}
 The area of unit vector field $V$ is expressed as an integral depending on the latitude $\alpha$ and the derivatives of the angle function $\theta$, as given by equation \eqref{Eq: Vol_derivatives}: 
\begin{eqnarray}\label{Eq:Area_Derivatives}
    \area\left(V\right) = \int\limits_{A_{\alpha_0}}{\sqrt{1+  \big( \tan \alpha + \theta_1\big)^2 + (\theta_2)^2}\, \, \,\nu},
\end{eqnarray}
where $\nu$ denotes the area form in the annulus region $A_{\alpha_0}$.

\noindent According to hypothesis \textit{(iii)} in Definition \ref{Def:AreaMinimizing_Annulus}, the equation \eqref{Eq:Area_Derivatives} takes the form
\begin{eqnarray} \nonumber
\area\left(V\right)
&=&\int\limits_{A_{\alpha_{0}}}\sqrt{\sec^{2}(\alpha) +\left( \theta _{2}\right) ^{2}}\,\, \nu  \\ \nonumber
&=& \int\limits_{-\alpha_0}^{\alpha_0} \int\limits_{0}^{2\pi} \cos (\alpha) \sqrt{\sec ^{2}(\alpha) +\left( \theta
_{2}\right) ^{2}}\, \, d\beta d\alpha  \\ 
&=& \int\limits_{-\alpha_0}^{\alpha_0} \int\limits_{0}^{2\pi} \sqrt{1+\cos^{2}(\alpha) \left( \theta _{2}\right) ^{2}}\, \, d\beta d\alpha . \label{Eq:Area(V)}
\end{eqnarray}

\noindent Taking $J := 1 + \cos^2( \alpha) \, (\theta_2)^2.$ By straightforward computation one obtains $J= H + I$, where
\begin{eqnarray}\nonumber
    H &=& \left( \sqrt{1 - \cos^2 (\alpha_0) \sec^2 (\alpha)} + \cos (\alpha_0) \, (\theta_2) \right)^2, \\
    I &=& \left( -\cos (\alpha_0) \sec (\alpha) + \sqrt{1 - \cos^2 (\alpha_0) \sec^2 (\alpha)} \, \cos (\alpha) \, (\theta_2) \right)^2. \label{Eq:Definicao I}
\end{eqnarray}
Therefore, the following inequality holds
\begin{eqnarray} \nonumber
\area\left(V\right)
&\geq& \int\limits_{-\alpha_0}^{\alpha_0} \int\limits_{0}^{2\pi} \sqrt{1 - \cos (\alpha_0) \sec ^{2}(\alpha)}\, \, d\beta d\alpha  + \cos (\alpha_0) \int\limits_{-\alpha_0}^{\alpha_0} \int\limits_{0}^{2\pi} \left( \theta _{2}\right) ^{2}\, \, d\beta d\alpha \\
&=& K_{\alpha_0} + \cos (\alpha_0) \, 2\pi^2, 
\end{eqnarray}
which concludes the proof of \textit{(i)}.

\noindent For all \(\alpha \in \left[-\alpha_0, \alpha_0\right]\), observe that for the equality in \textit{(i)} to hold, it is required that \(I\), as defined in \eqref{Eq:Definicao I}, equals zero, or equivalently
\begin{eqnarray*}
    \theta_2 = \dfrac{\cos (\alpha_0) \sec (\alpha)}{\sqrt{1 - \cos^2 (\alpha_0) \sec^2 (\alpha)} \, \cos (\alpha)}.
\end{eqnarray*}

\noindent Since \(\theta_2\) is the derivative of the angle function in the direction of the tangent vector to the meridian, we can write
\begin{eqnarray*}
    d\theta_p(e_2) = \theta_2(p) = \dfrac{d\theta (p)}{d\alpha},
\end{eqnarray*}
for all $p \in A_{\alpha_0}$. It follows that
\begin{eqnarray*}
    \theta(p) & = & \cos (\alpha_0)\int\dfrac{ \sec (\alpha)}{\sqrt{1 - \cos^2 (\alpha_0) \sec^2 (\alpha)} \, \cos (\alpha)} \,\, d\alpha\\
    & = & \arcsin\bigg(\cot(\alpha_0)\tan(\alpha)\bigg) + C,
\end{eqnarray*}
where $C$ is a real constant. For all $\beta \in (0,2\pi]$, let us define the condition
\begin{equation*}
\theta \left( 0,\beta \right) :=\frac{\pi }{2}.
\end{equation*}
Thus, 
\begin{equation*}
\theta \left( 0,\beta \right) = \arcsin(0)+C=\frac{\pi }{2}
\end{equation*}
which implies that $C=\dfrac{\pi }{2}.$ Therefore,   
\begin{equation*}
\theta \left( \alpha ,\beta \right) = \arcsin\bigg(\cot(\alpha_0)\tan(\alpha)\bigg) +\frac{\pi }{2}.
\end{equation*}
A direct verification shows that
\begin{equation*}
\theta \left( \alpha _{0},\beta \right)  = \pi   \quad \mbox{and}  \quad 
\theta \left( -\alpha _{0},\beta \right)  = 0,
\end{equation*}
which concludes the proof of \textit{(ii)}.

\end{proof}

\newpage



\section*{References}

\begin{itemize}
\item[[$\,\, 1$]]
Vincent Borrelli and Olga Gil-Medrano.
\newblock Area-minimizing vector fields on round 2-spheres.
\newblock \emph{J. Reine Angew. Math.}, 640:85–99, 2010.

\item[[$\,\, 2$]]
Fabiano G. B. Brito, Pablo M. Chacón, and David L. Johnson.
\newblock Unit vector fields on antipodally punctured spheres: big index, big volume.
\newblock \emph{Bull. Soc. Math. France}, 136(1):147–157, 2008.

\item[[$\,\, 3$]]
Fabiano G. B. Brito, Jackeline Conrado, Icaro Gonçalves, and Adriana V. Nicoli.
\newblock Area minimizing unit vector fields on antipodally punctured unit 2-sphere.
\newblock \emph{Comptes Rendus Mathématique}, 359-10:1225–1232, 2021.

\item[[$\,\, 4$]]
Fabiano G. B. Brito, André O. Gomes, and Icaro Gonçalves.
\newblock Poincaré index and the volume functional of unit vector fields on punctured spheres.
\newblock \emph{Manuscripta Math.}, 161(3-4):487–499, 2019.

\item[[$\,\, 5$]]
Sorin Dragomir and Domenico Perrone.
\newblock \emph{Harmonic vector fields}.
\newblock Elsevier, Inc., Amsterdam, 2012.

\item[[$\,\, 6$]]
Herman Gluck and Wolfgang Ziller.
\newblock On the volume of a unit vector field on the three-sphere.
\newblock \emph{Comment. Math. Helv.}, 61(2):177–192, 1986.

\item[[$\,\, 7$]]
Shigeo Sasaki.
\newblock On the differential geometry of tangent bundles of Riemannian manifolds.
\newblock \emph{Tohoku Math. J. (2)}, 10:338–354, 1958.
\end{itemize}

\end{document}